\documentclass[12pt]{article}
\usepackage[active]{srcltx}
\usepackage{amsmath}
\usepackage{amssymb}
\usepackage[usenames, dvipsnames]{color}
\newtheorem{theorem}{Theorem}
\newtheorem{corollary}{Corollary}

\usepackage{graphicx}
\DeclareGraphicsExtensions{eps,tif,emf}
\begin{document}
\begin{center}
\textbf{A fast algorithm for solving a three-dimensional inverse multiple frequency problems of scalar acoustics in a cylindrical region}
\end{center}
\begin{center}
{\bf Anatoly B. Bakushinsky}\footnote{Federal Research Center Computer Science and Control of Russian Academy of Sciences, Institute for Systems Analysis, 117312 Moscow; Mari State University, Lenin Sqr. 1, 424000 Yoshkar-Ola}, {\bf Alexander S. Leonov}\footnote{Department of Mathematics, National Nuclear Research University 'MEPHI',
Kashirskoe Shosse 31, 115409 Moscow}
\end{center}

%


\begin{abstract}
A new algorithm for the stable solution of a three-dimensional scalar inverse problem of acoustic sounding of an inhomogeneous medium in a cylindrical region is proposed. The data of the problem is the complex amplitude of the wave field, measured outside the region of acoustic inhomogeneities in a cylindrical layer. Using the Fourier transform and Fourier series, the inverse problem is reduced to solving a set of one-dimensional Fredholm integral equations of the first kind, to the subsequent calculation of the complex amplitude of the wave field in the region of inhomogeneity, and then to finding the required sound velocity field in this region. The algorithm allows solving the inverse problem on a personal computer of average performance for sufficiently fine three-dimensional grids in tens of seconds. A numerical study of the accuracy of the proposed algorithm for solving model inverse problems at various frequencies is carried out, and the issues of stability of the algorithm with respect to data perturbations are investigated.
\end{abstract}

\textbf{Keywords:} Three-dimensional wave equation, inverse
coefficient problem, regularizing algorithm, fast
Fourier transform.

\textbf{AMS Mathematics Subject Classification:} 65R20, 65R30, 65R32








\section{Introduction}\label{s1}
Let the scalar function $p(\textbf{x},t)$ define an acoustic wave field depending on the coordinates $\textbf{x}=(x,y,z)$ and time $t\ge 0$ in a domain $Q\subset {\bf {\mathbb R}}^{3} $. The domain is an infinite cylinder of the form $Q=\{(x,y,z): \, x^2+y^2\le b^2,\,z\in \mathbb{R}\}$.
The field is created by sources localized in a known region $S$. The medium is characterized by the local phase velocity of sound $c(\textbf{x})$ and has a constant density. Moreover, it is known that $c(\textbf{x})=c_{0} =\mathrm{const}$ outside the given region $X, X\subset Q$, $X\cap S=\varnothing $. In the region $ X $, the function $c(\textbf{x})$ can vary, and this is interpreted as the location of acoustic inhomogeneities there. In this case, the field $p(\textbf{x},t)$ for a harmonic source of the form $f(\textbf{x},\omega )e^{i\omega t} $ with a known frequency $\omega $ can be found within the linear acoustics approximation as $p(\textbf{x},t)=u(\textbf{x},\omega )e^{i\omega t} $, where the complex amplitude $u(\textbf{x},\omega )$ satisfies the equation
\begin{equation}\label{L1_1_}
\Delta u(\textbf{x},\omega )+k_{0}^2 u(\textbf{x},\omega )=f(\textbf{x},\omega )+\omega ^{2}
\xi (\textbf{x})u(\textbf{x},\omega ),\, \, \, \textbf{x}\in Q .  
\end{equation}
Here $k_{0} =\frac{\omega }{c_{0} } $, and $\xi (\textbf{x})=c_{0}^{-2} -c^{-2} (\textbf{x})$.
We also assume that the boundary condition
\begin{equation*}
\left(\frac{
\partial u}{\partial \textbf{n}} +\sigma(x)u\right)_{\partial Q} =0
\end{equation*}
is satisfied with known function $\sigma(\textbf{x})$, 
and the radiation condition is fulfilled in the coordinate $z$. Here $\textbf{n}$ is the outer normal to the surface $\partial Q$. Without going into details of the conditions for the coefficients $\sigma(\textbf{x})$ and $f(\textbf{x},\omega )$, we make the following assumption.

\textbf{Assumption 1.} The function $\xi (\textbf{x})$ is continuous with compact support in $ X $ and the corresponding problem \eqref{L1_1_} with the indicated additional conditions has a unique solution $u(\textbf{x},\omega )\in H_1(Q)$ for each considered $\omega$.

Finding such a function $u(\textbf{x},\omega )$ is
\emph{direct problem}. We are interested in the following \emph{inverse problem} for the equation  with the indicated additional conditions: knowing the complex amplitude of the field $u(\textbf{x},\omega )$ for some set of frequencies $\omega $ in the domain $ Y$ ($Y\subset Q,\,Y\cap X=\varnothing, \,Y\cap S=\varnothing  $) find the coefficient $\xi (\textbf{x})$, i.e. the function $c(\textbf{x})$ defining acoustic inhomogeneities in the region $ X $. Introducing the Green's function $G(\textbf{x},\textbf{x}',\omega )$ for the Helmholtz equation in the domain $ Q $, we can reduce the inverse problem under certain assumptions about the smoothness of the functions $u(\textbf{x},\omega ),f(\textbf{x},\omega ),\, c(\textbf{x})$ to a nonlinear system of integral equations for the unknowns $u(\textbf{x}',\omega )$ and $\xi(\textbf{x}'),\,\textbf{x}'\in X$:
\begin{eqnarray}\label{L1_2_}
&&u({\bf{x}},\omega ) = {u_0}({\bf{x}},\omega ) + {\omega ^2}\int_X^{} {G({\bf{x}} , {\bf{x}}',\omega )\xi ({\bf{x}}')u({\bf{x}}',\omega )d{\bf{x}}'} ,\,\,\,{\bf{x}} \in X,\nonumber\\&&  \\
&&{\omega ^2}\int_X^{} {G({\bf{x}} , {\bf{x}}',\omega )\xi ({\bf{x}}')u({\bf{x}}',\omega )d{\bf{x}}' = w({\bf{x}},\omega ),\,\,{\bf{x}} \in Y}\nonumber
\end{eqnarray}
(see, for example, \cite{22,21,1}, etc.). The functions
\[u_{0} (\textbf{x},\omega )=\int _{X}^{}G(\textbf{x},\textbf{x}',\omega )f(\textbf{x}',\omega )d\textbf{x}' ,\, \, \, w(\textbf{x},\omega )=u(\textbf{x},\omega )-u_{0} (\textbf{x},\omega ),\, \, \textbf{x}\in Y,\]
included in \eqref{L1_2_} are known (computable) functions, and the values
$u(\textbf{x}',\omega )$ and $\xi (\textbf{x}'),\, \textbf{x}'\in X$
must be defined.

The problems similar to \eqref{L1_2_} has been sufficiently well investigated theoretically. In particular, questions of the existence and uniqueness of their solution for various domains $ Q, X, Y $ have been studied (see, for example, \cite{21}--\cite{2} and etc.). We emphasize that the purpose of this article is not to study the properties of the inverse problem \eqref{L1_2_}.
We only propose an effective numerical method for solving it in the case of cylindrical regions $ X, Y $. That is why we do not carry out a detailed reduction of the inverse problem to the system  \eqref{L1_2_} indicating all the requirements for the coefficients, but make only one more assumption.

\textbf{Assumption 2.} Reduction of the inverse problem to the \eqref{L1_2_} system is possible for given $u_0(\textbf{x},\omega )\in L_2(Q)$. The system \eqref{L1_2_} is solvable and defines, for all considered $ \omega $, some solution of the inverse problem, that is a continuous function $\xi (\textbf{x})$ with compact support in $ X $ and the function $u(\textbf{x},\omega )\in L_2(X)$.

Note that under such assumption the system \eqref{L1_2_} may have more than one solution.

In many papers, various numerical methods for solving problems similar to \eqref{L1_2_} in two-dimensional and three-dimensional formulations are consi\-de\-red. For example, in the work \cite{1} the system \eqref{L1_2_} is reduced to a nonlinear operator equation that is then solved by a special iterative method. This approach does not take into account explicitly that this equation is an ill-posed problem, and the iterative method used is not regularized. Nevertheless, the method works for model scatterers of medium strength \cite[p.101-105]{1}. In the work \cite{2}, the regularized Gauss -- Newton method was used to solve the system of equations \eqref{L1_2_} in the three-dimensional axially symmetric case. A special gradient method and the Fletcher-Reeves method were used for a similar problem in the paper \cite{3}. Other gradient techniques have been used in \cite{4},\cite{5}.

Note that there are also alternative approaches to solving the inverse problem of acoustic sounding that are not related to a system of the type \eqref{L1_2_}. In particular, the original \emph{method of boundary control} was proposed and developed in the works \cite{6},\cite{7} that allows solving three-dimensional problems \eqref{L1_2_}. In \cite{8}, the well-known method of R.G. Novikov \cite{9} was applied to solve the two-dimensional inverse acoustic scattering problem. In the subsequent work \cite{10}, a comparative analysis of a variant of this method and some other functional-analytical methods for solving two-dimensional inverse problems of acoustic scattering was carried out. Methods by M.V. Klibanov summarized in the monograph \cite{11}, as well as methods from the monograph \cite{Kab}, turned out to be very promising in processing real experimental data. Also of great interest are recent works \cite{Kl1},\cite{Kl2}.

All the above methods for solving the inverse problem of acoustic sounding require significant computational resources in a three-dimensional formulation. In this regard, we note the articles \cite{12},\cite{13}, in which the original inverse coefficient problem for the wave equation, from which the equation \eqref{L1_1_} is actually obtained, is reduced to a three-dimensional Fredholm integral equation of the first kind with the right-hand side containing special integrals of the recorded field. The proposed method for solving this equation turned out to be very effective numerically, and allows one to quickly solve three-dimensional inverse problems for sufficiently fine grids on a personal computer (PC) even without parallelization. This method is also very useful for the \eqref{L1_2_} problem, and we will demonstrate it below.

In this article, we adhere to the following scheme for solving a nonlinear system\eqref{L1_2_} for each considered $\omega$.

1) Introducing the notation $v(\textbf{x}',\omega )=\xi (\textbf{x}')u(\textbf{x}',\omega ),\, \, \textbf{x}'\in X$, solve the second equation (linear integral Fredholm equation of the first kind), written in the form
\begin{equation} \label{L1_3_} \omega ^{2} \int _{X}^{}G(\textbf{x},\textbf{x}',\omega )v(\textbf{x}',\omega )d\textbf{x}'=w(\textbf{x},\omega ),\, \, \textbf{x}\in Y,
\end{equation}
for the functions $v(\textbf{x}',\omega )$.

 2) With found function $v(\textbf{x}',\omega )$, calculate the function $u(\textbf{x},\omega ),\, \, \textbf{x}\in X,$ from the first equality of the system \eqref{L1_2_}, written in the form
\begin{equation} \label{L1_4_}
u(\textbf{x},\omega )=u_{0} (\textbf{x},\omega )+\omega ^{2} \int _{X}^{}G(\textbf{x},\textbf{x}',\omega )v(\textbf{x}',\omega )d\textbf{x}' ,\, \, \, \textbf{x}\in X.
\end{equation}

 3) Find the function $\xi (\textbf{x})$ from the equation $v(\textbf{x},\omega )=\xi (\textbf{x})u(\textbf{x},\omega ),\, \, \textbf{x}\in X$ using computed functions $v(\textbf{x},\omega ),u(\textbf{x},\omega )$.

A similar scheme has been used before. For example, in \cite{Smi}, it was used to solve an ill-posed problem, a three-dimensional equation similar to \eqref{L1_3_}. However, regularization methods were not used there, and the solution of the inverse problem was sought for a rather narrow class of functions with a piecewise constant current.

Below we will show that under some not very burdensome special assumptions about the set of data registration $ Y $ and about the set $ X $, where the solution of the inverse problem is sought, scheme 1) - 3) is effectively implemented numerically and allows solving the corresponding three-dimensional inverse problem for fairly fine grids even on a personal computer of average performance in a few tens of seconds without parallelization. The proposed algorithm for solving the inverse problem and its numerical study are the main results of this work.

The article is organized as follows. In Sect.2, we consider the geometric scheme used and reduction of the inverse problem to a final system of equations. In Sect.3, a method for obtaining model data for the inverse problem (Algorithm 1) is given, and Algorithm 2 for solving our inverse problem in a cylindrical domain is presented and discussed. Sect.4 is devoted to a finite-dimensional approximation and solution of model inverse problems. Here we consider some details of algorithms and present results of numerical experiments. The properties of the algorithm, such as the speed of finding a solution and sensitivity of solutions to input data errors, are discussed in Sect.5. Finally, in Sect.6 we formulate main conclusions.

\section{Geometric scheme used and further reduction of the inverse problem}
Everywhere below, the domains $ X, Y $ have the form of cylinders: the solution domain is $X=\{ 0\le r\le a\} \times {\bf {\mathbb R}}_{z}$, $r=\sqrt{x^{2} +y^{2}} $, and the observation region is $Y=\{ r_{0} \le r\le b\} \times {\bf {\mathbb R}}_{z} $.
Figure 1 shows the disposition of these regions. The possible position of the sources is also conventionally shown there. We consider the cylinders to be infinite in the variable $z$: $z\in (-\infty ,+\infty
)$, bearing in mind that the sought-for function $ \xi $ is compactly supported and in particular with respect to this variable. In what follows we will denote as $X_{xy}, Y_{xy}$ sections of cylindrical areas $ X, Y $ by plane, perpendicular to the $Oz$ axis.
\begin{figure}[h]
  \centering
\includegraphics[width=120mm,height=90mm]{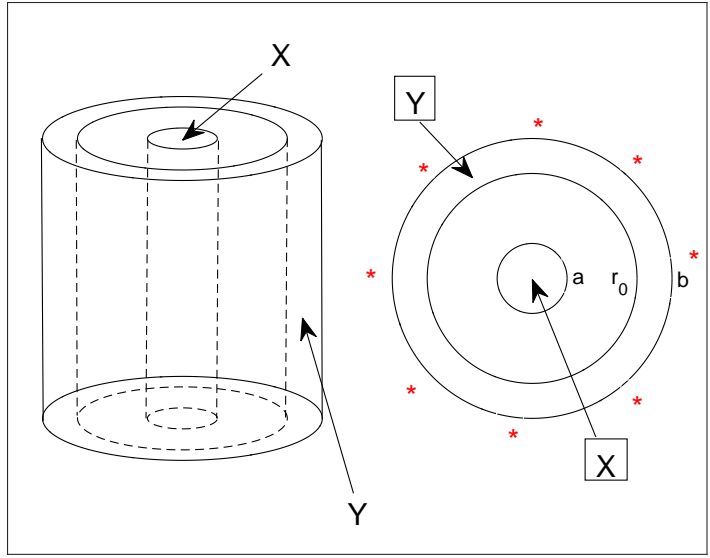}
  \caption{{\small Geometric scheme of data registration for the inverse problem: $ X $ is the domain of wave field scatterers, $ Y $ is the region of data registration, asterisks are the conditional positions of the sources.}}
  \label{fig6}
\end{figure}

Green's function for the problem \eqref{L1_1_} with the indicated additional conditions in the cylindrical domain $ Q $ is well known,
and here we do not write it out explicitly. Note an important feature of this function: in cylindrical coordinates $\mathbf{x}=(r,\varphi,z),\mathbf{x}'=(r',\varphi',z')$ it has the form:
\begin{multline*}
 G\left( { {{\bf{x}} , {\bf{x}}'} ,\omega } \right) = G\left( {\sqrt {{r^2} + r{'^2} - 2rr'\cos (\varphi  - \varphi ') + {{(z - z')}^2}} ;\omega } \right) = \\={G_0}\left( {r,r',\cos (\varphi  - \varphi '),z - z';\omega } \right),
\end{multline*}
and the expression for the function $ {G_0} $ can be found, for example, in
\cite{Budak}. In the same coordinates, we represent the functions
$u({\bf{x}},\omega ) = u(r,\varphi ,z),\,v({\bf{x}},\omega ) = v(r,\varphi ,z)$. Then the equations \eqref{L1_3_}, \eqref{L1_4_} can be written as follows:
\begin{multline}\label{eq21}
{\omega ^2}\int\limits_0^a {\int\limits_0^{2\pi } {\int\limits_{ - \infty }^{ + \infty } {{G_0}\left( {r,r',\cos (\varphi  - \varphi '),z - z';\omega } \right)} } } v(r',\varphi ',z';\omega)r'dr'd\varphi' dz'=\\=w(r,\varphi ,z;\omega),\,\,r\in[r_0,b], \varphi\in[0,2\pi],\,z\in \mathbb{R},
\end{multline}
and
\begin{multline}\label{eq21-1}
u(r,\varphi ,z;\omega) - {u_0}(r,\varphi ,z;\omega) = \\={\omega ^2}\int\limits_0^a {\int\limits_0^{2\pi } {\int\limits_{ - \infty }^{ + \infty } {{G_0}\left( {r,r',\cos (\varphi  - \varphi '),z - z';\omega } \right)} } } v(r',\varphi ',z';\omega)r'dr'd\varphi' dz',\\
r\in[0,a], \varphi\in[0,2\pi],\,z\in \mathbb{R}.
\end{multline}
Here $v(r',\varphi ',z';\omega)=\xi (r',\varphi ',z')u(r',\varphi' ,z';\omega)$. Next, we introduce the Fourier transforms $F_{z}[\cdot ](\Omega)=\int_{-\infty }^{+\infty }[\cdot ]e^{i\Omega z}dz$ with respect to the variable $ z $ (or $ z' $) for functions $G_0,v,u,u_0,w$ as elements of space $ L_{2} $ (see Assumption 2):
\[\begin{array}{l} {\tilde{G}(r,r',\cos \varphi,\Omega;\omega )=F_{z}
\left[{G_0}\left( {r,r',\cos \varphi,z;\omega } \right)\right](\Omega)},\\ {\tilde{w}(r,\varphi ,\Omega;\omega)=F_{z} \left[w(r,\varphi ,z;\omega)\right](\Omega )},
 \\ {\tilde{u}(r,\varphi ,\Omega;\omega)=F_{z} \left[u(r,\varphi ,z;\omega)\right]
(\Omega )},\, {\tilde{u}_0(r,\varphi ,\Omega;\omega)=F_{z} \left[u_0(r,\varphi ,z;\omega)\right]
(\Omega ).}\end{array}\]
Then, by the convolution theorem, the equality \eqref{eq21} and \eqref{eq21-1} can be written as
\begin{multline}\label{eq22}
{\omega ^2}\int\limits_0^a {\int\limits_0^{2\pi } {\tilde G(r,r',\cos (\varphi  - \varphi '),\Omega ;\omega )\tilde v(r',\,} } \varphi ',\Omega ;\omega )r'dr'd\varphi ' =\\= \tilde w(r,\varphi ,\Omega ;\omega ),\,r \in [{r_0},b],
\end{multline}
\begin{multline}\label{eq221}
\tilde u(r,\varphi ,\Omega ;\omega ) - {{\tilde u}_0}(r,\varphi ,\Omega ;\omega ) = \\={\omega ^2}\int\limits_0^a {\int\limits_0^{2\pi } {\tilde G(r,r',\cos (\varphi  - \varphi '),\Omega ;\omega )\tilde v(r',\,} } \varphi ',\Omega ;\omega )r'dr'd\varphi ',\,
r \in [0,a],
\end{multline}
with $\tilde{v}(r',\varphi ',\Omega;\omega)=F_{z} \left[v(r',\varphi ',z';\omega)\right](\Omega )=F_{z} \left[\xi (r',\varphi ',z')u(r',\varphi' ,z';\omega)\right]$. We also introduce expansions of the functions $\tilde G,\,\tilde v,\,\tilde u,\,{{\tilde u}_0},\tilde w$ in the basis system of functions $\{ e^{in\varphi } \} ,\, n\in {\bf {\mathbb Z}},$ in the space $L_{2} (0,2\pi )$:
\[\begin{array}{l}
\tilde G(r,r',\cos \varphi ,\Omega ;\omega ) = \sum\limits_n {{G_n}(r,r',\Omega ;\omega )} {e^{in\varphi }},\\\tilde v(r',\varphi ',\Omega ;\omega ) = \sum\limits_n {{v_n}(r',\Omega ;\omega )} {e^{in\varphi '}},\\
\tilde u(r,\varphi ,\Omega ;\omega ) = \sum\limits_n {{u_n}(r',\Omega ;\omega )} {e^{in\varphi }},\,{{\tilde u}_0}(r,\varphi ,\Omega ;\omega ) = \sum\limits_n {{u_{0n}}(r',\Omega ;\omega )} {e^{in\varphi }};\,\,\\
\tilde w(r,\varphi ,\Omega ;\omega ) = \sum\limits_n {{w_n}(r',\Omega ;\omega )} {e^{in\varphi }};\,\,\varphi ,\varphi ' \in [0,2\pi ]
\end{array}\]
with coefficients
\[\begin{array}{l}
{G_n}(r,r',\Omega ;\omega ) = \frac{1}{{2\pi }}\int_0^{2\pi } {\tilde G(r,r',\cos \varphi ,\Omega ;\omega ){e^{ - in\varphi }}d\varphi ,\,} \,\\
{v_n}(r',\Omega ;\omega ) = \frac{1}{{2\pi }}\int_0^{2\pi } {\tilde v(r',\varphi ',\Omega ;\omega ){e^{ - in\varphi '}}d\varphi '},\\{{w_n}(r,\Omega ;\omega ) = \frac{1}{{2\pi }}\int_0^{2\pi } {\tilde w(r,\varphi ,\Omega ;\omega ){e^{ - in\varphi }}d\varphi ,\,} } \\
{u_n}(r,\Omega ;\omega ) = \frac{1}{{2\pi }}\int_0^{2\pi } {\tilde u(r,\varphi ,\Omega ;\omega ){e^{ - in\varphi }}d\varphi ,} \\{u_{0n}}(r,\Omega ;\omega ) = \frac{1}{{2\pi }}\int_0^{2\pi } {{{\tilde u}_0}(r,\varphi ,\Omega ;\omega ){e^{ - in\varphi }}d\varphi. \,}
\end{array}\]
Then the relations \eqref{eq22}, \eqref{eq221} can be reduced to the following system of equalities that is true for any $n \in\mathbb{ Z},\, \Omega  \in \mathbb{R}$ and all considered $\omega$:
\begin{equation}\label{eq23}
{\omega ^2}\int_0^a {{G_n}(r,r',\Omega ;\omega ){v_n}(r',\Omega ;\omega )r'dr'}  = \frac{1}{{2\pi }}{w_n}(r,\Omega ;\omega ),\,r \in [{r_0},b],
\end{equation}
\begin{multline}\label{eq24}
\frac{1}{{2\pi }}\left( {{u_n}(r,\Omega ;\omega ) -
{u_{0n}}(r,\Omega ;\omega )} \right) = {\omega ^2}\int_0^a
{G_n}(r,r',\Omega ;\omega ){v_n}(r',\Omega ;\omega )r'dr',\\r \in [0,a].
\end{multline}
Moreover, the equality
\begin{multline}\label{eq25}
 {v_n}(r',\Omega ;\omega ) =  \frac{1}{{2\pi }}\int_0^{2\pi } {{F_z}\left[ {v(r',\varphi ,z;\omega )} \right]{e^{ - in\varphi }}d\varphi \,}  = \\
 = \frac{1}{{2\pi }}\int_0^{2\pi } {{F_z}\left[ {\xi (r',\varphi ,z)u(r',\varphi ,z;\omega )} \right]{e^{ - in\varphi }}d\varphi \,}  =\\= \frac{1}{{2\pi }}\int_0^{2\pi } {{F_z}\left[ {\xi (r',\varphi ,z)\sum\limits_m {{u_m}(r',\Omega ;\omega )} {e^{im\varphi }}} \right]{e^{ - in\varphi }}d\varphi}
\end{multline}
is satisfied. The relations \eqref{eq23}, \eqref{eq24} are equations for the unknown functions ${v_n}(r',\Omega ;\omega )$ and ${u_n}(r,\Omega ;\omega )$, depending on one variable $ r'$ or $ r $. Other arguments of these functions, i.e. $ \Omega, \omega $, are parameters.

\section{Algorithms for solving direct and inverse problems in a cylindrical domain}
\subsection{Direct problem}
The direct problem considered below is to find the function ${w}\left(r,\varphi,z,\omega\right),\, r\in [r_0,b],\varphi\in [0,2\pi],z\in\mathbb{ R}$, from the equalities \eqref{eq21}, \eqref{eq21-1}
using the known finite function $\xi (r',\varphi',z')$ and the given function of sources ${u_0}(r,\varphi ,z;\omega)$ for a set of frequencies $ \omega $ under consideration. For this, the equalities \eqref{eq21}, \eqref{eq21-1} are reduced to the system \eqref{eq22},\eqref{eq221}, and then to the system of relations \eqref{eq23} -- \eqref{eq25}. As a result, the calculation of the function $ w $ can be represented as the following algorithm.

\textbf{Algorithm 1}

Step 1). For the set of frequencies $ \omega $ under consideration, calculate the Fourier transforms in $ z $:
\begin{multline*}
  \tilde{G}(r,r',\cos \varphi,\Omega;\omega )=F_{z}
\left[{G_0}\left( {r,r',\cos \varphi,z;\omega } \right)\right](\Omega ), \\
  \tilde{u}_0(r,\varphi ,\Omega;\omega)=F_{z} \left[u_0(r,\varphi ,z;\omega)\right](\Omega )
\end{multline*}
and expand the resulting functions into Fourier series in the variable $\varphi$:
\begin{multline*}
\tilde G(r,r',\cos \varphi ,\Omega ;\omega ) = \sum\limits_n {{G_n}(r,r',\Omega ;\omega )} {e^{in\varphi }},\\{{\tilde u}_0}(r,\varphi ,\Omega ;\omega ) = \sum\limits_n {{u_{0n}}(r,r',\Omega ;\omega )} {e^{in\varphi }}.
\end{multline*}
Both of these procedures can be implemented using the Fast Discrete Fourier Transform (FFT).

Step 2). For each parameters $\omega ,\Omega$, we implement the following iterative process of solving the equations \eqref{eq24}, \eqref{eq25} with respect to the set of function $\{{u_n}(r,\Omega ;\omega )\}$, $n\in \mathbb{Z}$:
\begin{multline}\label{L1_6_}
 v_n^{(k)}(r',\Omega ;\omega ) =  \frac{1}{{2\pi }}\int_0^{2\pi } {{F_z}\left[ {\xi (r',\varphi ,z)\sum\limits_m {u_m^{(k)}(r',\Omega ;\omega )} {e^{im\varphi }}} \right]{e^{ - in\varphi }}d\varphi},\\r' \in [0,a],
\end{multline}
\begin{multline}\label{L1_66}
 {u_n^{(k+1)}(r,\Omega ;\omega ) ={u_{0n}}(r,\Omega ;\omega )}+ 2\pi{\omega ^2}\int_0^a {G_n}(r,r',\Omega ;\omega )v_n^{(k)}(r',\Omega ;\omega )r'dr',\\r \in [0,a],\,\,k=0,1,2,...,
\end{multline}
with an initial guess $\{u_n^{(0)}(r,\Omega ;\omega )\}=\{{u_{0n}}(r,\Omega ;\omega )\}$.

Step 3). Stop the process by some rule at iteration number $ \nu $ and obtain an approximate solution $\{u_n^{(\nu)}(r,\Omega ;\omega )\}$ of the system of equations \eqref{eq24}, \eqref{eq25}, and related functions $\{v_n^{(\nu)}(r,\Omega ;\omega )\}$.

Step 4). Calculate an approximate set of values $\{{w_n}(r,\Omega ;\omega )\}$ from \eqref{eq23} :
\begin{equation}\label{L1_60}
 w_n^{(\nu)}(r,\Omega ;\omega )=2\pi {\omega ^2}\int_0^a {{G_n}(r,r',\Omega ;\omega )v_n^{(\nu)}(r',\Omega ;\omega )r'dr'},\,r \in [{r_0},b]
\end{equation}
and accept the function
\begin{equation}\label{L1_77}
\tilde w^{(\nu)}(r,\varphi ,\Omega ;\omega ) = \sum\limits_n {w_n^{(\nu)}(r',\Omega ;\omega )} {e^{in\varphi }}
\end{equation}
or its inverse Fourier transform as an approximate solution of the direct problem.

We will not carry out a theoretical analysis of the convergence of Algorithm 1 here, since we do not formally use it in solving the inverse problem. It is only needed to generate the model data $w$. We only note that, as follows from the general theory of solving integral equations of the second kind (see, for example, \cite{15}), the algorithm will rapidly converge, at least for small $ \omega $. Below we will demonstrate numerical examples confirming this statement.

\subsection{Inverse problem}
The task is to find for each frequency $\omega$ the solution $\xi (r',\varphi',z')$ of the system \eqref{eq21}, \eqref{eq21-1} using given function ${w}\left(r,\varphi,z,\omega\right),\, r\in [r_0,b],\varphi\in [0,2\pi],z\in\mathbb{ R}$, and given source function ${u_0}(r,\varphi ,z;\omega)$. To do this, we reduce the problem \eqref{eq21}, \eqref{eq21-1} to the system \eqref{eq23} -- \eqref{eq25}, assuming the sets of data functions $\{{w_n}(r,\Omega ;\omega )\}$ and $\{{u_{0n}}(r,\Omega ;\omega )\}$ are calculated. Then from the system and these sets we find the function $ \xi $. The solution procedure is presented in the form of the following algorithm.

\textbf{Algorithm 2}

Step 1). For all considered parameters $\omega ,\Omega $ and all used $ n $ we solve one-dimensional integral equations of the first kind corresponding to the equalities \eqref{eq23}:
\begin{equation}\label{L1_8_}
{\omega ^2}\int_0^a {{G_n}(r,r',\Omega ;\omega ){v_n}(r',\Omega ;\omega )r'dr'}  = \frac{1}{{2\pi }}{w_n}(r,\Omega ;\omega ),\,r \in [{r_0},b].
\end{equation}
Here we use a suitable regularization method (regularizing algorithm, RA) for these ill-posed problems. The result is a set of approximate solutions $\{{v_n}(r',\Omega ;\omega )\}$.

Step 2). Using the found functions $\{{v_n}(r',\Omega ;\omega )\}$, we calculate from the equalities \eqref{eq24} the set of functions
$\{u_n(r,\Omega ;\omega )\}$:
\begin{equation}\label{L1_9_}
{u_n(r,\Omega ;\omega ) ={u_{0n}}(r,\Omega ;\omega )}+ 2\pi{\omega ^2}\int_0^a {{G_n}(r,r',\Omega ;\omega )v_n(r',\Omega ;\omega )r'dr',\,r \in [0,a].}
\end{equation}

Step 3). We restore the functions $v(r,\varphi,z,\omega )$ and $u(r,\varphi,z,\omega )$ for $r\in [0,a],\varphi\in [0,2\pi],z\in \mathbb{R}$ using the sets $\{{v_n}(r,\Omega ;\omega )\}$ and $\{u_n(r,\Omega ;\omega )\}$, summing the corresponding Fourier series in the variable $\varphi$ and then calculating inverse Fourier transform in the variable $z$, $F_\Omega^{-1}[\cdot]$:
\begin{multline*} v(r,\varphi ,z;\omega)=F_{\Omega }^{-1} \left[\sum _{n}v_{n} (r,\Omega;\omega )e^{in\varphi }  \right](z),\\
 u(r,\varphi ,z;\omega)=F_{\Omega }^{-1} \left[\sum _{n}u_{n} (r,\Omega;\omega )e^{in\varphi }  \right](z),\, (r,\varphi ,z)\in X,
\end{multline*}

Step 4). Next, we find the solution $\xi (r,\varphi,z)$ from the equation \[u(r,\varphi,z,\omega )\xi (r,\varphi,z)=v(r,\varphi,z,\omega )\] for each point $(r,\varphi,z)\in X$. This can be done for each $\omega $ under consideration, and the result will generally depend on $\omega $. Further, for example, one can average the results by one or another method over the value of $\omega $.

Step 5). Finally, we calculate the function $c(r,\varphi,z)$ from the equality $\xi (r,\varphi,z)=c_{0}^{-2} -c^{-2} (r,\varphi,z)$.

Let's make some comments on Algorithm 2.

a) To implement Step 1 it is necessary to clarify the properties of the functions ${v},\,{w}$. Using Assumption 2, we presume that the inclusions ${v}(\textbf{x},\omega)\in L_2(X)$, ${w}(\textbf{x},\omega)\in L_2(Y)$ are valid for each considered frequency $ \omega $.
The first inclusion follows from the compactness of the support of the function $ \xi $, and the second is postulated. In this case, the known methods for solving linear ill-posed problems in Hilbert spaces are applicable to the equations \eqref{L1_8_} (see, for example, \cite{2},\cite{3},\cite{16},\cite{17},\cite{18}, etc.).

{b) In solving our inverse problem, the equation \eqref{L1_3_}, i.e. \eqref{eq21} and the equations \eqref{L1_8_} generated by it, may have more than one solution for the finite set of frequencies $ \omega $ used. Therefore, it is important to establish a connection between the solutions of all these equations. The connection is substantiated by the following statements.

\begin{theorem}\label{T1}
1) Let $w(r,\varphi ,z;\omega)\in L_{2} (Y)$ for every $ \omega $. Then any solution $v(r,\varphi ,z;\omega)\in L_{2} (X)$ to the equation \eqref{eq21} can be represented in the form
\begin{equation} \label{eq__7_} v(r,\varphi ,z;\omega)=F_{\Omega }^{-1} \left[\sum _{n}v_{n} (r,\Omega;\omega )e^{in\varphi }  \right](z),\, \, (r,\varphi ,z)\in X,
\end{equation}
where the functions $v_{n} (r,\Omega;\omega )\in L_{2} \left\{[0,a]\times
{\bf {\mathbb R}}_{\Omega } \right\}$ satisfy the integral equations \eqref{eq23} for each $\omega$. Conversely, if $v_{n}
(r,\Omega;\omega )\in L_{2} \left\{[0,a]\times {\bf {\mathbb R}}_{\Omega
} \right\}$ are solutions of equations \eqref{eq23} such that $\sum
_{n}\left\| v_{n} (r,\Omega;\omega )\right\| _{L_{2}
\left\{[0,a]\times {\bf {\mathbb R}}_{\Omega } \right\}}^{2}
<\infty $ for each $ \omega $, then a function of the form \eqref{eq__7_} is the solution to the equation \eqref{eq21}.

2) For each $ \omega $ the equality
\begin{equation}\label{eq__6_}
\left\| v(r,\varphi ,z;\omega)\right\| _{L_{2} (X)}^{2} =
\sum _{n}\left\| v_{n} (r,\Omega;\omega)\right\| _{L_{2} \left\{[0,a]\times {\bf {\mathbb R}}_{\Omega } \right\}}^{2}
\end{equation}
holds.
\end{theorem}

\begin{corollary}\label{C1}
Let the functions $\bar{v}_{n} (r,\Omega;\omega )\in L_{2}
\left\{[0,a]\times {\bf {\mathbb R}}_{\Omega } \right\}$ be normal solutions to equations \eqref{eq23} (solutions with minimal norm). Then the function
\[\bar{v}(r,\varphi ,z;\omega)=F_{\Omega }^{-1} \left[\sum
_{n}\bar{v}_{n} (r,\Omega;\omega )e^{in\varphi }  \right](z),\, \,
(r,\varphi ,z)\in X,\] is the unique normal solution to the equation
\eqref{eq21}.
\end{corollary}

The proofs of these statements are carried out in the same way as in \cite{12}. There, similar statements were proved for solutions of another integral equation of the first kind, similar to \eqref{eq21} in form and properties and differing only in the kernel and the right-hand side. For brevity, we omit repeating these proofs.

c) In the used scheme for solving the inverse problem and in Algorithm 2, the inversion procedure is actually used only when solving equations of the first kind \eqref{eq23} or \eqref{L1_8_}, which are ill-posed problems.
Due to the possible non-uniqueness of solutions to these equations, we apply methods aimed at finding normal solutions (Tikhonov regularization, TSVD method) to implement Step 1 of the algorithm. If the solution is unique, it coincides with the calculated normal solution. The indicated regularization methods have been substantiated and tested in a number of works (see for example, \cite{2},\cite{3},\cite{16},\cite{17},\cite{18} and etc.). For similar equations, such a regularization was used in \cite{12},\cite{13}.
Note that Step 1 is the most laborious step when using Algorithm 2.

\section{Finite-dimensional approximation and solution of model problems}
Everywhere below, it is assumed that the equations \eqref{L1_3_}, \eqref{L1_4_} and their consequences are written in dimensionless form with $c_{0} =1$, so that $k_{0} =\omega $. Model domains are as follows:
\begin{multline*}
\hspace{-3mm}Q=\left\{ (x,y,z):~x^{2}+y^{2}\leq 4^{2},|z|\leq 2\right\}, X = \left\{ {(x,y,z):\,{x^2} + {y^2} \leq1,\,|z| \leqslant 2} \right\}, \\
 Y = \left\{ {(x,y,z):\,{3^2} \leq {x^2} + {y^2} \leqslant {4^2},\,|z| \leq2} \right\}.
\end{multline*}
It is also supposed that $\sigma(\textbf{x})=0$ and model sources are given in the form \[f(\textbf{x})=\sum _{m=1}^{M_{s}}A_{m} \delta (\textbf{x}-\textbf{x}_{m} ),\] where $\textbf{x}_{m} $ are coordinates of $ \delta $-shaped point sources. Then \[u_{0} (\textbf{x},\omega )=\sum _{m=1}^{M_{s}}A_{m} G_0(\textbf{x},\textbf{x}'_{m},\omega ),\] and the Fourier transform of this function, $\tilde{u}_{0} (r,\varphi,\Omega;\omega  )$, can be calculated in advance.
We did not set ourselves the goal of optimizing the number, positions and amplitudes of sources. In all calculations, it was assumed that ${M_{s}} = 8,\,{A_m} = 1,\,{{\mathbf{x}}_m} = ({r_m},{\varphi _m},{z_m})$ with ${r_m} = 4.01$ and
\[
 {\varphi _m} = \left[ {0,\,\frac{\pi }{2}, - \frac{\pi }{2},\pi ,0,\,\frac{\pi }{2}, - \frac{\pi }{2},\pi } \right], \,{z_m} = [ - 1, - 1, - 1, - 1,1,1,1,1].
\]
For the numerical study of the proposed algorithms, two direct and inverse model problems are considered. \emph{The first problem} has a solution $\xi (x,y,z)$  of the form
 \begin{multline*}
   \xi (x,y,z)=\xi _{1}(x,y,z)+\xi _{2}(x,y,z),~(x,y,z)\in Q; \\
   \xi _{1}(x,y,z)=\left\{ A_{0}\exp \left\{ -30R_{1}(x,y,z)\right\} ,~(x,y,z)\in Q_{1};~0,~(x,y,z)\in Q\backslash Q_{1}\right\},\\
   \xi _{2}(x,y,z)=\left\{ 2A_{0}\exp \left\{ -30R_{2}(x,y,z)\right\} ,~(x,y,z)\in Q_{2};~0,~(x,y,z)\in Q\backslash Q_{2}\right\},
 \end{multline*}
with
\begin{multline*}
  R_{1}(x,y,z)=5(x-0.4)^{2}+5y^{2}+0.125(z+0.1)^{2},\\
  R_{2}(x,y,z)=5(x+0.4)^{2}+5(y-0.4)^{2}+0.125(z-0.2)^{2}
\end{multline*}
and
\begin{multline*}
  Q_{1}=\left\{ (x,y,z):~(x-0.4)^{2}+y^{2}+0.125(z+0.1)^{2}\leq 1.3^{2}\right\} \\
  Q_{2}=\left\{ (x,y,z):~(x+0.4)^{2}+(y-0.4)^{2}+0.125(z-0.2)^{2}\leq 0.5^{2}\right\}
\end{multline*}
This function simulates small local inhomogeneities of the medium, the position of which and the corresponding velocity distributions must be found. Algorithm 2 is tuned specifically to search for such inhomogeneities.
The value $A_{0} $ determines the contrast \[\frac{\Delta c}{c_{0} } =\underset{\textbf{x}}{\max } \left\{\frac{1}{\sqrt{1-c_{0}^{2} \xi (\textbf{x})} } \right\}-1\] of the desired solution. In the calculations, we used $A_{0} =0.545$, and this corresponds to a contrast equal to 15.95. According to the classification from \cite[c.33]{1}, such a scatterer can be considered strong, taking into account its characteristic dimensions $l\sim 0.2$ (see Fig.\ref{fig4}A) and the values $c_0=1,\omega=3$: $\frac{\Delta c}{c_{0} }\gg \frac{c_0}{l\omega}$.
Such scatterers are quite common in practice. The second model task will be described below.

The equations \eqref{eq21}, \eqref{eq21-1} were approximated in the domains $ X, Y $ by the finite-difference method on uniform grids of
variables $r,r',\varphi ,z$. The sizes of the grids are determined by the numbers $N_r,N_{r'},N_\varphi,N_z$. In the region $ X $, the grid has size $N_{r'}\times N_\varphi\times N_z$, and in the region $ Y$, the size is $N_{r}\times N_\varphi\times N_z$. Specific dimensions will be given below for each example. Discrete analogues of the functions $\tilde G(r,r',\cos \varphi ,\Omega ;\omega )$, ${{\tilde u}_0}(r,\varphi ,\Omega ;\omega )$, used in Algorithms 1 and 2, were calculated for the considered frequencies $ \omega $ from the known values ${G_0}\left( {r,r',\cos \varphi ,z;\omega } \right)$, ${u_0}(r,\varphi ,z;\omega )$ using the fast Fourier transform with grid $\{\Omega^{(m)} \}_{m=1}^{N_z}$ in the variable $ \Omega $. Fourier series expansions were also implemented using the FFT with $n\in [0,N_\varphi-1]$. Details of these well-known standard calculations are given, for example, in \cite{16}.

\subsection{Obtaining model data for the inverse problem using Algorithm 1}
Now we present typical results of a numerical study of the iterative process \eqref{L1_6_},\eqref{L1_66} to obtain data for the first inverse problem on grids of size $N_r=32,N_{r'}=33,N_\varphi=90,N_z=64$. Figure 2 shows a comparison of the convergence rate for the process with various quantities
$\omega =k_{0} $.
\begin{figure}
  \centering
\includegraphics[width=80mm,height=60mm]{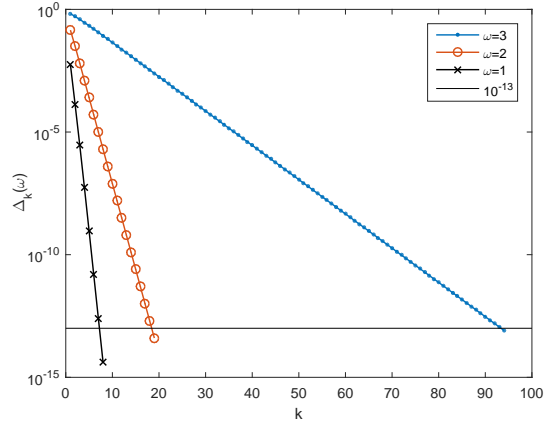}
  \caption{{\small Convergence rate of iterations \eqref{L1_6_}, \eqref{L1_66} for different $\omega =k_{0} $.}}
  \label{fig2}
\end{figure}

\noindent The iterations were stopped by the value \[\Delta_k(\omega)=\frac{{\left\{ {\sum\limits_n {\mathop {\left\| {u_n^{(k )}\left( {r,\Omega ;\omega } \right) - u_n^{(k -1)}\left( {r,\Omega ;\omega } \right)} \right\|}\nolimits_{{L_2}\left\{ {\Pi} \right\}}^2 } } \right\}^{1/2}}}{{\left\{ {\sum\limits_n {\mathop {\left\| {u_n^{(0)}\left( {r,\Omega ;\omega } \right)} \right\|}\nolimits_{{L_2}\left\{ {\Pi} \right\}}^2 } } \right\}^{1/2}}}
,\]
when the condition
$\Delta_\nu(\omega) \leqslant {10^{ - 13}}$ was satisfied for the iteration number $ \nu $. Here $\Pi=\{(r,\Omega)\in[0,a] \times {\mathbb{R}_\Omega }\}$. After that, using the found set of functions $\left\{ {u_n^{(\nu )}\left( {r,\Omega ;\omega } \right)} \right\}$ the functions $\left\{ {v_n^{(\nu )}(r',\Omega ;\omega )} \right\}$ were calculated by the formula \eqref{L1_6_} with $ k = \nu $. Then, using the formula \eqref{L1_60}, the functions $\left\{ {w_n^{(\nu )}(r,\Omega ;\omega )} \right\}$ were found. Further, they were transformed according to \eqref{L1_77}, in the function ${{\tilde w}^{(\nu )}}(r,\varphi ,\Omega ;\omega )$. The inverse Fourier transform of the last function with respect to the variable $ z $, ${w^{(\nu )}}(r',\varphi ,z;\omega )$, represents the data for solving the inverse problem.
The form of this function, found for $\omega =k_{0} =3$, is shown in Fig.\ref{fig3} for $ z = 0 $.
\begin{figure}
  \centering
\includegraphics[width=130mm,height=70mm]{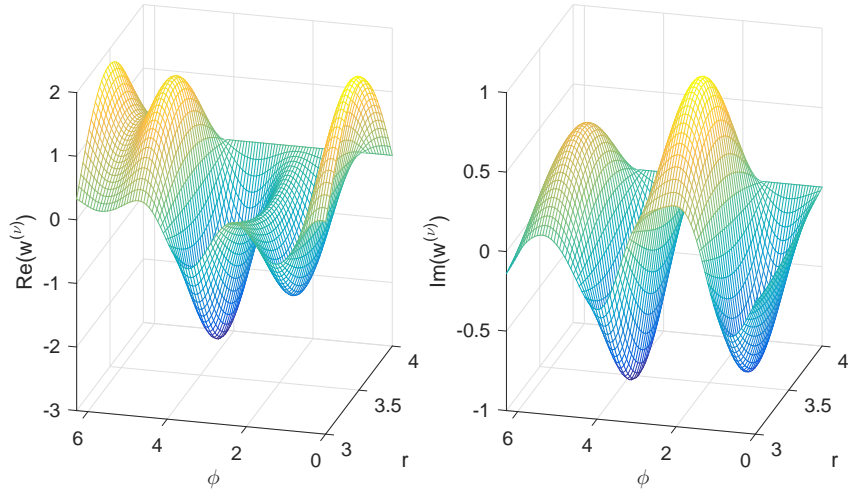}
  \caption{{\small Typical data $w^{(\nu)}(r,\varphi,0;\omega =3)$ for solving the inverse problem.}}
  \label{fig3}
\end{figure}

Further, the data for solving the inverse problem were specified with some perturbations that can be interpreted as a measurement error. In our calculations, this was modeled by imposing an additive normally distributed pseudo-random noise with zero mean on the function $w^{(\nu )}(r,\varphi,z;\omega )$ so that the resulting approximate function $w^{(\nu )}_{\delta } (r,\varphi,z;\omega )$ would satisfy the condition \[\left\| w^{(\nu) }_{\delta } (r,\varphi,z;\omega )-w^{(\nu) } (r,\varphi,z;\omega )\right\| _{L_{2} (Y)} \le \delta \left\| w^{(\nu )} (r,\varphi,z;\omega )\right\| _{L_{2} (Y)}. \] This corresponds to approximate data with relative error $\delta $.

\subsection{Implementation of Algorithm 2}
The first step of Algorithm 2, i.e. solving equations of the first kind \eqref{L1_8_} by regularization methods was discussed in the papers \cite{12},\cite{13} in connection with the solution of another inverse problem, formally similar to the considered one and differing from it only in the form of the kernel and the right-hand side. In these papers, it was noted that for each considered frequency $ \omega $, the discretization used reduces the equations \eqref{L1_8_} to a system of linear algebraic equations (SLAE) of the form $A^{(m)}_n {V}^{(m)}_n ={W}^{(m)}_n $ for each $n\in[0,N_\varphi-1]$, $m\in[1,N_z]$. Here $A^{(m)}_n =\left[\mu _{ij} {G_n} (r_i,r'_j,\Omega ^{(m)};\omega )\right]_{i=1, j=1}^{N_r, N_{r'}}$ is the matrix of the system obtained by discretizing the kernel of the equation \eqref{L1_8_} on the considered grid of size $N_r\times N_{r'} $, quantities $\Omega ^{(m)} $ are grid points along $ \Omega $, and $\mu _{ij} $ are quadrature coefficients for calculating integrals in \eqref{L1_8_}. The right sides of the system, ${W}^{(m)}_n =\frac{1}{2\pi\omega^2}\left[{w}_{n}\left(r_i,\Omega ^{(m)};\omega  \right)\right]_{i=1}^{N_r}$, are column vectors of height $ N_r $, and the column vector ${V}^{(m)}_n$ of height $ N_r' $ contains unknowns $v_n(r'_j,\Omega ^{(m)};\omega)$. Thus, when performing Step 1 of Algorithm 2 for each $ \omega $ it is necessary to solve $N_\varphi\times N_z$ systems of linear equations with matrices of size $N_r\times N_{r'}$.

We solved the indicated SLAEs using various versions of the Tikhonov regularization method \cite{16},\cite{17} and using the TSVD  method \cite{18}. The justi\-fi\-cation for these methods is given in \cite{12},\cite{13} for similar problem. The best calculation results were obtained using the TSVD method. We present them further.

Step 2 of Algorithm 2 does not cause any difficulties for the discretized problem, since is reduced to matrix multiplication of discrete quantities ${G}_n,{v}_n$ and addition of the result with a discrete analogue of the function ${u}_{0n} $. Step 3 was performed using the inverse FFT. Finally, Step 4 was implemented for each considered frequency $ \omega $ using the following procedure for finding the normal pseudosolution of the equation $u\xi=v$ by the TSVD method at each point $(r',\varphi',z)\in X$: $\xi =\left\{ \frac{v}{u},~|u|>\mathrm{tol};0,~|u|\leq \mathrm{tol}\right\} $ with $\mathrm{tol}=10^{-12}$. Further, it is easy to recalculate the function $\xi(\mathbf{x})$ into $c(\mathbf{x})$. For the sake of brevity, we do not do this in the examples below, presenting the value $\xi(\mathbf{x})$ directly in the figures.
\begin{figure}
  \centering
\includegraphics[width=140mm,height=85mm]{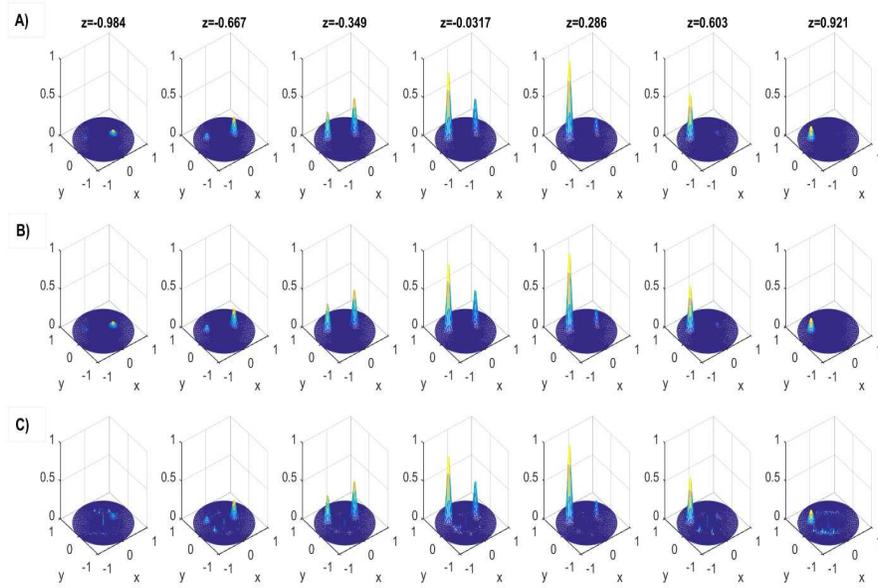}
  \caption{{\small Qualitative comparison of the exact solution $\xi _{\rm exact} (x,y,z)$ and the calculated approximate solutions $\xi _{\rm appr} (x,y,z)$ of the inverse problem in different sections $z=\mathrm{const}$.
   A) exact solution; B) approximate solution obtained for exact data by Algorithm 2; C) approximate solution for disturbed data with $\delta =10^{-8} $.}}
  \label{fig4}
\end{figure}
\begin{figure}
  \centering
\includegraphics[width=90mm,height=70mm]{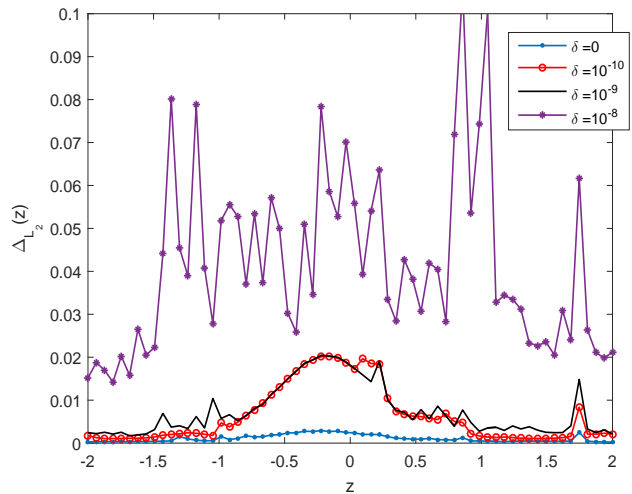}
  \caption{{\small The first model problem for $\omega=3$. Relative error $ \Delta_{L_ {2}}(z) $ of approximate solutions for different $ z $ at different levels of data perturbation $ \delta $.  }}
  \label{fig5}
\end{figure}
\begin{figure}
  \centering
\includegraphics[width=60mm,height=50mm]{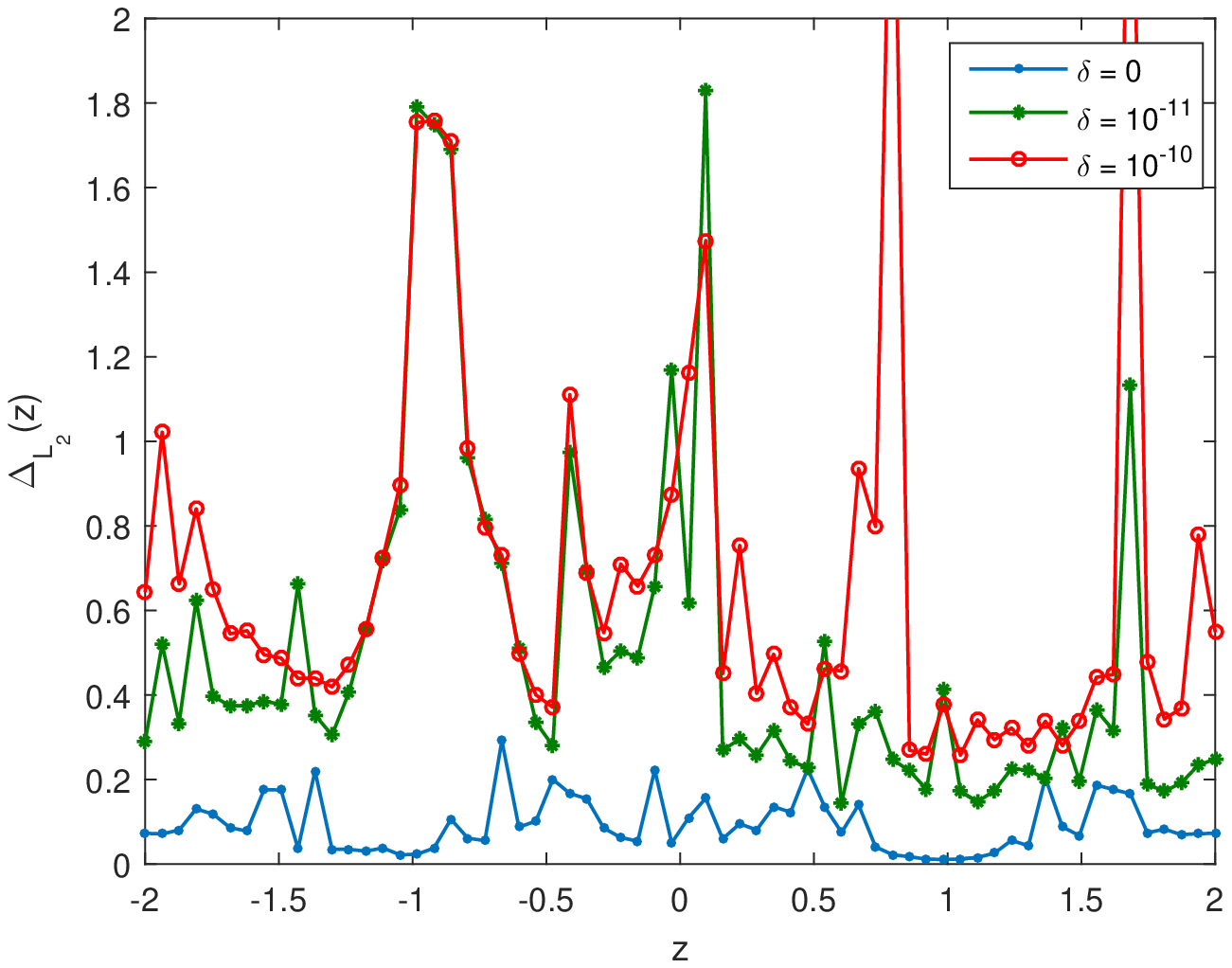}
\includegraphics[width=60mm,height=50mm]{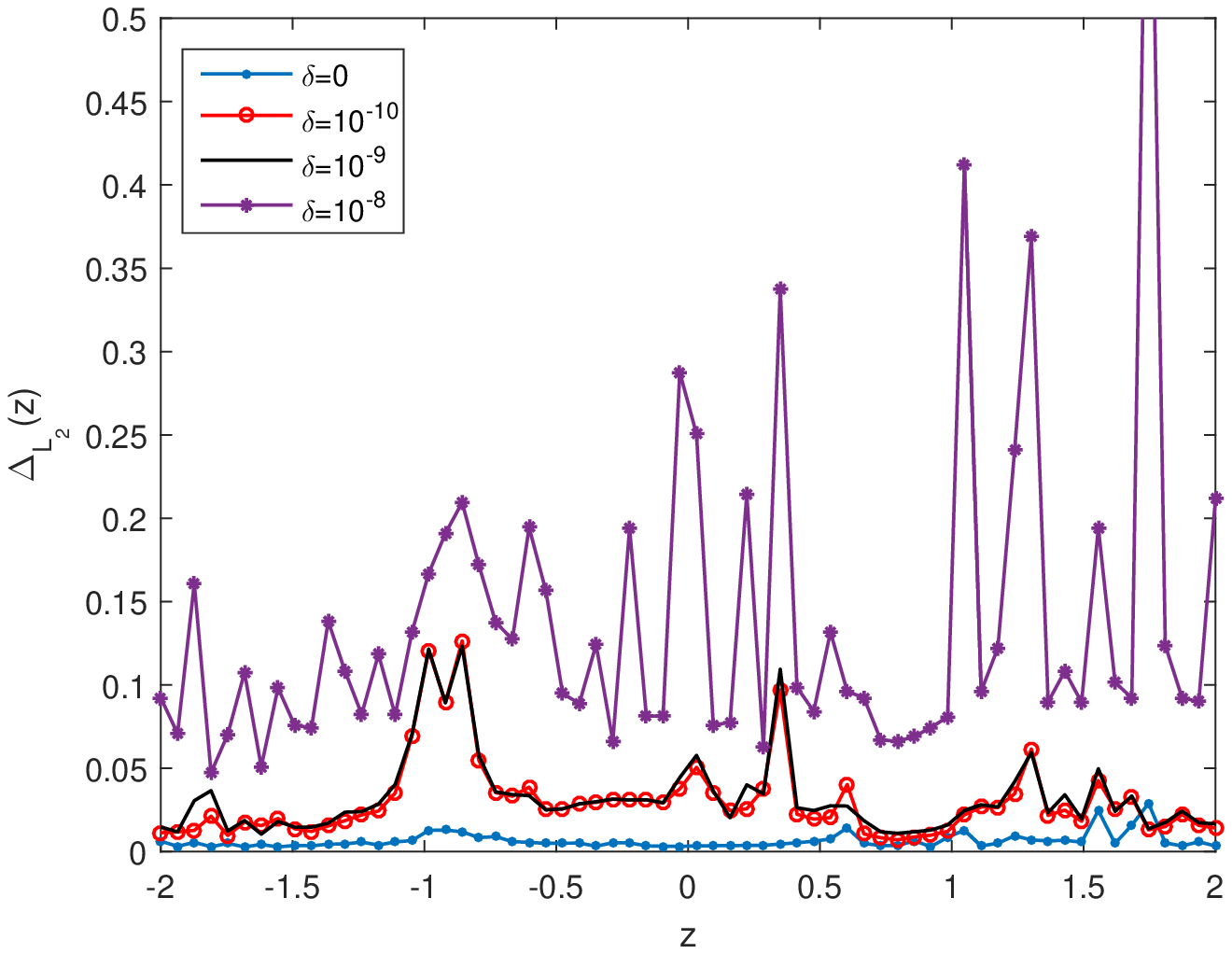}
  \caption{{\small The first model problem. Relative error $\Delta _{L_{2} } (z)$ of approximate solutions for different $ z $ at different levels of data perturbation $ \delta $. Left: for $ \omega = 1 $; right: for $ \omega = 2 $.}}
  \label{fig6}
\end{figure}
\begin{figure}
  \centering
\includegraphics[width=110mm,height=80mm]{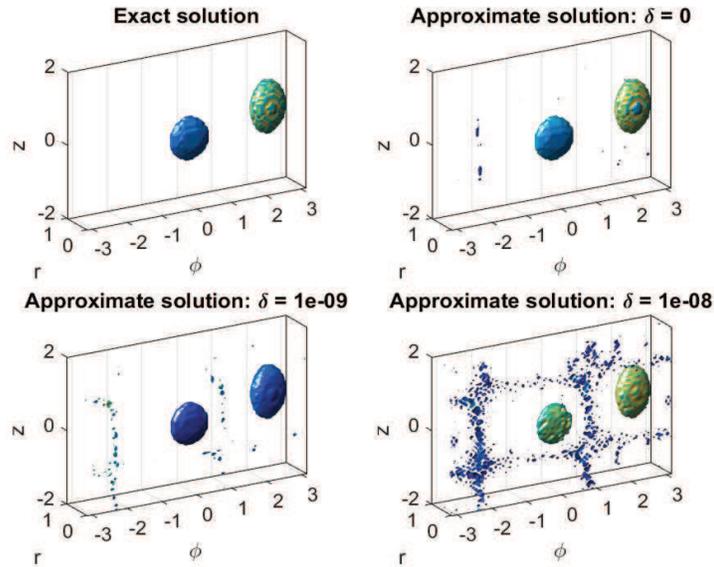}
  \caption{{\small Second model problem. A qualitative comparison of the positions and geometry of the reconstructed inhomogeneity $\xi(\mathbf{x})$ for different $\delta$.}}
  \label{fig7}
\end{figure}
\begin{figure}
  \centering
\includegraphics[width=100mm,height=80mm]{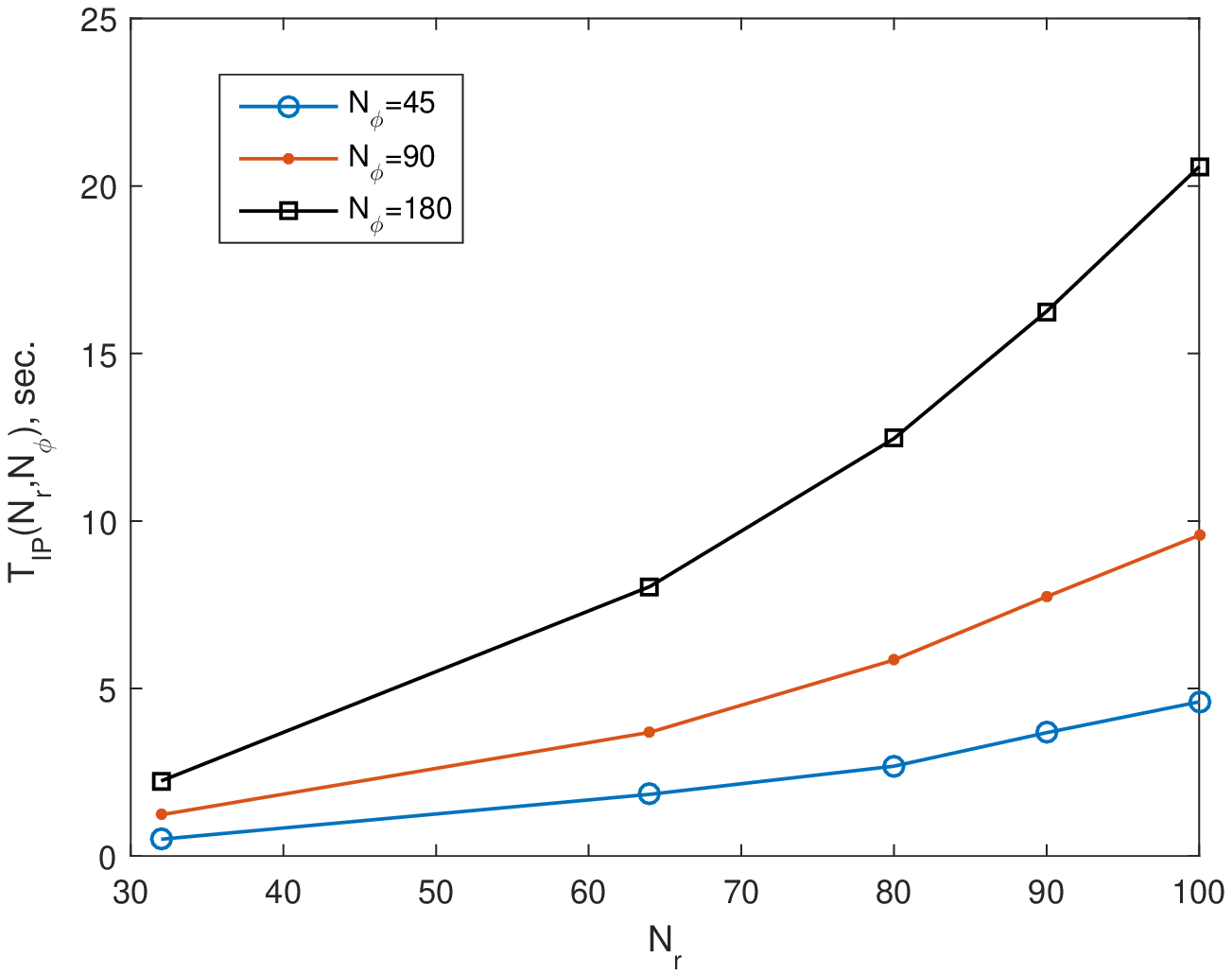}
  \caption{{\small Time $T_{IP} (N_r,N_\varphi)$ of solving the inverse problem for different $N_r,N_\varphi$.}}
  \label{fig8}
\end{figure}

\subsection{Results of solving the inverse problem}
We start by applying Algorithm 2 in numerical experiments of solving the first model inverse problem on grids of size $N_r=32,N_{r'}=33,N_\varphi=90,N_z=64$ for $\omega=k_0=3$.
The problem was solved with exact data, more precisely, calculated by Algorithm 1, and approximate data with different error levels $\delta $. For a qualitative comparison, Fig.\ref{fig4} shows the exact and approximate solutions of the inverse problem, $\xi _{\rm exact} (x,y,z)$ and $\xi _{\rm appr} (x,y,z)$, in different sections with $ z =\mathrm{const}$.
The exact solution is shown in the first line of the figure. The second line contains approximate solution for exact data. The third line shows the solution for the perturbed data with $\delta =10^{-8} $. The figure demonstrates a fairly high sensitivity of solutions to data disturbances. More detailed information on the accuracy of solving the inverse problem, i.e. on the relative error
\[\Delta _{L_{2} } (z)=\frac{\left\| \xi _{\rm appr} (x,y,z)-\xi _{\rm exact} (x,y,z)\right\| _{L_{2} (X_{xy})} }{\underset{z}{\max }\left\| \xi _{\rm exact} (x,y,z)\right\| _{L_{2} (X_{xy})} } \]
of the approximate solution $\xi _{\rm appr} (x,y,z)$ in the sections $z=\mathrm{const}$, is presented in Fig.\ref{fig5} for different levels $ \delta $ of data error.

For comparison, Fig.\ref{fig6} shows the errors of approximate solutions obtained for different $ \delta $ in solving model problems with $ \omega = 1 $ and $ \omega = 2 $. The improvement in accuracy is obvious with increasing number $\omega$.

When solving the inverse problem under consideration, it is very important to know how accurately the algorithm allows one to determine the positions of the investigated local scatterers. For illustration, \emph{the second model problem} was solved, which differs from the first one only in expressions for $\xi _{1}(x,y,z)$ and $\xi _{2}(x,y,z)$:
\begin{multline*}
  \xi _{1}(x,y,z)=\left\{ A_{0};~(x,y,z)\in Q_{1};~0,~(x,y,z)\in Q\backslash Q_{1}\right\}\\
 \xi _{2}(x,y,z)=\left\{ 2A_{0};~(x,y,z)\in Q_{2};~0,~(x,y,z)\in Q\backslash Q_{2}\right\},\,\,A_0=0.2.
\end{multline*}
Such a solution corresponds to two ellipsoidal scatterers lying in the $ X $ region and filled with substances having different refractive indices.
Fig.\ref{fig7} shows qualitatively the influence of the problem data perturbation on the determination of the position and geometry of inhomogeneities. The contrast of the exact solution is 0.291. It is seen that the positions can be determined quite accurately in the case of $\delta =10^{-9} $. This is also possible for $\delta =10^{-8} $ when using appropriate noise filtering in the found solution.

\section{Some properties of Algorithm 2}

All calculations were carried out in MATLAB on a PC with an Intel (R) Core (TM) i7-7700 CPU 3.60 GHz, 16GB RAM without parallelization. Algorithm 2 for solving the inverse problem turned out to be fast enough. We present the results of corresponding numerical experiments for solving the first inverse problem with $\omega=2 $. In the experiments, it was assumed that the $ z $ grid is fixed ($N_z=64$), and only the sizes $N_r,N_\varphi$ of the $ r $ and $ \varphi $ grids change. In addition, it was assumed that $N_{r'}=N_r+1$.
Then the time for solving the inverse problem is a function of the form $T_{IP} (N_r,N_\varphi)$. This dependence is shown in Fig.\ref{fig8}. When the grid size in the variable $ z $ is changed, the time $T_{IP} (N_r,N_\varphi)$ changes proportionally to the number $ N_z $, since the time is determined by solving $ N_z \times N_\varphi $ equations of the form \eqref{L1_8_}.

Note once again that the inverse problem being solved is very sensitive to input data errors. When solving it with double precision, introducing random errors with an amplitude of the order of $ 10^{-8} $ into the right-hand sides of Eqs.\eqref{L1_8_} leads to serious distortions of the solution. This happens when using the TSVD method and the regularization method too. This sensitivity is associated with a very fast decrease in the singular numbers of the matrices $A^{(m)}_n $ for SLAEs solved in Step 1 of Algorithm 2, and this is a specific feature of the inverse problem being solved.
 A similar property of the inverse coefficient problem for the wave equation was noted before in the works \cite{12,13}. The corresponding theoretical estimates of the error under various a priori assumptions on the exact solution can be found in \cite{2,3}.

\section{Conclusions}

From the numerical experiments carried out in this work, the following conclusions can be drawn.

1. The considered three-dimensional inverse problem of scalar acoustics in a cylindrical region can be solved using Algorithm 2 for sufficiently fine grids in a few tens of seconds on a PC of average performance even without parallelization. For this, one should use the scheme for recording data of the inverse problem in a cylindrical layer indicated in the article. The proposed algorithm can be easily parallelized.

2. The inverse problem under consideration is in itself very sensitive to data perturbations; to obtain a detailed approximate solution, data measured with high accuracy are required. This feature of the problem does not depend on the used algorithm.

3. Algorithm 2 makes it possible to reliably determine the position and shape of small local inhomogeneities of the acoustic medium with data having small errors.\vspace{3mm}

{\bf Acknowledgments}

 This work was supported by the Russian Science Foundation (project 20-11-20085) for the first author in part of substantiating numerical algorithms and the Programm of Competitiveness Increase of
the National Research Nuclear University MEPhI (Moscow Engineering
Physics Institute); contract no. 02.a03.21.0005, 27.08.2013 for the second author.


\begin{thebibliography}{2}
\bibitem{22} Ramm A.G. \emph{Multidimensional Inverse Scattering Problems,} Pitman Monogr. Surv. Pure Appl. Math. 51. Harlow: Longman Scientific \& Technical, 1992.

\bibitem{21} Colton D., Kress R. \emph{Inverse Acoustic and Electromagnetic Scattering Theory,} 2nd ed.. Appl. Math. Sci. 93. Berlin: Springer, 1998.

\bibitem{1} Goryunov A.A., Saskovets A.V. \emph{Inverse scattering problems in acoustics,} M., Publishing House of Moscow State University, 1989 (in Russian).

\bibitem{2} Bakushinsky A., Goncharsky A. \emph{Ill-Posed Problems: Theory and Applications.} Dordrecht: Kluwer Academic Publishers, 1994.

\bibitem{3} Bakushinsky A.B., Kokurin M.Yu. \emph{Iterative methods for approximate solution of inverse problems,} Mathematics and Its Applications. Dordrecht: Kluwer Academic Publishers, 2004.

\bibitem{4} A. V. Goncharsky and S. Y. Romanov, On two approaches to the solution of coefficient inverse problems for wave equations, Zh. Vychisl. Mat. Mat. Fiz. 52 (2012), no.2, pp.263-269.

\bibitem{5} A. V. Goncharsky and S. Y. Romanov, Supercomputer technologies in inverse problems of ultrasound tomography, Inverse Problems 29 (2013), no.7, Article ID 075004.

\bibitem{6} Belishev M.I. Recent progress in the boundary control method, Inverse Problems. 2007. V.23. N5. P.1-67.

\bibitem{7} Pestov L.N., Bolgova V.M., Danilin A.N. Numerical reconstruction of the threedimensional speed of sound by the method of boundary control, Bulletin of Ugra State University, 2011, Issue 3, pp.92-98 (in Russian).

\bibitem{8} Burov V.A., Alekseenko N.V., Rumyantseva O.D. Multifrequency Generalization of the Novikov Algorithm for the Two-Dimensional Inverse Scattering Problem, Acoustic Journal, V.55, No.6, 2009, pp.784-798.

\bibitem{9} Novikov P. G. Reconstruction of the two-dimensional Schrodinger operator from the scattering amplitude at a fixed energy, Funktsional. analysis and its adj., T.20, No.3, 1986, pp.90-91.

\bibitem{10} Burov V.A., Vecherin S.N., Morozov S.A., Rumyantseva O.D. Modeling of the Exact Solution of the Inverse Scattering Problem by Functional Methods. Acoustic Journal, Vol.56, No.4, 2010, pp.516-536.

\bibitem{11} Beilina L., Klibanov M.V., \emph{Approximate Global Convergence and
Adaptivity for Coefficient Inverse Problems,} New York: Springer, 2012. 

\bibitem{Kab} Kabanikhin S.I., Satybaev A.D., Shishlenin M.A. \emph{Direct Methods of Solving Multidimensional Inverse Hyperbolic Problems,} Utrecht: VSP,
    2004.

\bibitem{Kl1} Klibanov M.V., Kolesov A.E. Convexification of a 3-D coefficient inverse scattering problem // Computers and Mathematics with Applications. 2019. V.77. P.1681-1702.

\bibitem{Kl2} Klibanov M.V., Kolesov A.E., Nguyen Dinh-Liem. Convexification method for an inverse scattering problem and its performance for experimental backscatter data for buried targets, SIAM J. Imaging Sciences. 2019. V.12, N.1. P.576-603.

\bibitem{12} Bakushinsky A.B., Leonov A.S. Fast numerical method of solving 3D coefficient inverse problem for wave equation with integral data, Journal of Inverse and Ill-Posed Problems. 2018. V.26. Issue 4. P.477-492.

\bibitem{13} A.B. Bakushinskii, and A.S. Leonov, Low-Cost Numerical Method for Solving a Coefficient Inverse Problem for the Wave Equation in Three-Dimensional Space. Comp. Math. and Math. Phys., 2018, Vol. 58, No. 4, pp.548-561.

\bibitem{Smi} Evstigneev R.O., Medvedik M.Yu., Smirnov Yu.G., Tsupak A.A. The inverse problem of body's heterogeneity recovery for early diagnostics of diseases using microwave tomography, University proceedings, Volga region, Physical and Mathematical Sciences, 2017, No.4 (44), pp. 3-17 (in Russian).

\bibitem{Budak} B.M. Budak, A.A. Samarskii and A.N. Tikhonov, \emph{A Collection of Problems on Mathematical Physics,} Pergamon Press, Oxford, 1964


\bibitem{15} F. Riesz and B. Sz.-Nagy, \emph{Functional Analysis,} Frederick Ungar Publishing Co., N.Y., 1955.

\bibitem{16} A. N. Tikhonov, A. V. Goncharsky, V. V. Stepanov and A. G. Yagola, \emph{Numerical Methods for the Solution of Ill-Posed Problems}, Math. Appl. 328, Kluwer Academic Publishers, Dordrecht, 1995.

\bibitem{17} A. S. Leonov, \emph{Solution of Ill-Posed Inverse Problems. Theory Review, Practical Algorithms and MATLAB Demonstrations}, Librokom, Moscow, 2010, 2013 (in Russian).

\bibitem{18} Engl H.W., Hanke M, Neubauer A. \emph{Regularization of Inverse Problems}, Dordrecht: Kluwer Academic Publishers, 1996.


\end{thebibliography}
\end{document}